\documentclass{rspublic}
\pdfoutput=1
\usepackage[colorlinks=true]{hyperref}
\usepackage{graphicx}
\usepackage{amsmath}
\usepackage{amssymb}
\usepackage{color}

\newcommand{\bm}[1]{\boldsymbol{#1}}

\begin{document}

\title[Minimal surfaces bounded by elastic lines]{Minimal surfaces bounded by elastic lines}

\author[L. Giomi and L. Mahadevan]{L. Giomi and L. Mahadevan}

\affiliation{School of Engineering and Applied Sciences, and Department of Physics, Harvard University, Pierce Hall 29 Oxford Street Cambridge, MA 02138, USA.\\[7pt]} 

\label{firstpage}

\maketitle
 
\begin{abstract}{Soap films; minimal surfaces, Elastica}
In mathematics, the classical Plateau problem consists of finding the surface of least area that spans a given rigid boundary curve. A physical realization of the problem is obtained by dipping a stiff wire frame of some given shape in soapy water and then removing it; the shape of the spanning soap film is a solution to the Plateau problem. But what happens if a soap film spans a loop of inextensible but flexible wire? We consider this simple query that couples Plateau's problem to  Euler's {\em Elastica}: a special class of twist-free curves of given length that minimize their total squared curvature energy. The natural marriage of two of the oldest geometrical problems linking physics and mathematics leads to a quest for the shape of a minimal surface bounded by an elastic line: the Euler-Plateau problem. We use a combination of simple physical experiments with soap films that span soft filaments, scaling concepts, exact and asymptotic analysis combined with numerical simulations  to explore some of the richness of the shapes that result. Our study raises  questions of intrinsic interest in geometry and its natural links to a range of disciplines including materials science, polymer physics, architecture and even art.
\end{abstract}

\begin{center}
	\nointerlineskip
	\rule{0.9\textwidth}{0.7pt}
\end{center}

\section{Introduction}

Soap films are ethereal and ephemeral. Yet their beauty is not and has fascinated the scientist and the mathematician for more than two centuries, since Lagrange and Plateau first linked soap films and minimal surfaces: special surfaces of zero mean curvature that minimize the area enclosed by a given contour. The classical problem of finding the surface with least area that spans a given rigid boundary curve was first formulated by Lagrange (1760). However, the subject was lifted to a whole new level by Plateau in the mid nineteenth century using a series of beautiful experiments that showed that a physical realization of these objects naturally arises through a consideration of soap films (Plateau 1849). Since then, the subject has inspired mathematicians (Osserman 2002, Colding \& Minicozzi 2007, Morgan 2008), scientists (Thomas {\em et al}. 1988, Kamien 2002), engineers  and even artists (Ferguson) and spurred extraordinary developments in all these fields. Over the past two decades there has been a renewed appreciation for the subject as links to problems as disparate as materials science, string theory (Maldacena 1998) and other areas of high energy physics and cosmology (Penrose 1973) have become increasingly apparent.

Almost without exception, the focus of problems over the last two centuries has always been on determining the shape of a soap film with a given rigid boundary. However, since experiments with soap films are often carried out by dipping closed wire frames into soap and then pulling them out, there is a natural generalization of this boundary value problem that is suggested by the following question: what if the soap film is not bounded by a rigid contour, but has instead a soft boundary such as a flexible, inextensible wire? This simple question links at the hip the Plateau problem to another classical problem in geometry, that of the {\em Elastica}, a curve of given length that minimizes the total squared curvature, first formulated by Euler (1844) and since then object of continuous mathematical investigation (Bryant \& Griffiths 1983, Langer \& Singer 1984, Mumford 1993). The natural marriage of two of the oldest geometrical problems linking physics and mathematics leads to a quest for the shape of a minimal surface bounded by an elastic line. The competition between the surface tension of the film, that acts to minimize the area bounded by the loop, and the elasticity of the boundary, that tries to minimize its deformation, dictates the shape of this composite structure. Since small loops are stiff, they will remain flat and circular when spanned by a film. On increasing the size of the loop past some critical value, the forces due to surface tension can no longer be sustained by the boundary which buckles out of the plane even as the film takes on the shape of a saddle. Dimensional analysis suggests that the critical size of the loop when it just buckles scales as $(\alpha/\sigma)^{1/3}$ , where $\alpha=EI$ ($E$ Young's modulus and $I$ area moment of inertia of the elastic filament) is the bending stiffness of  the elastic line and $\sigma$ is the surface tension of the soap film.

\begin{figure}[t]
\centering
\includegraphics[width=1\columnwidth]{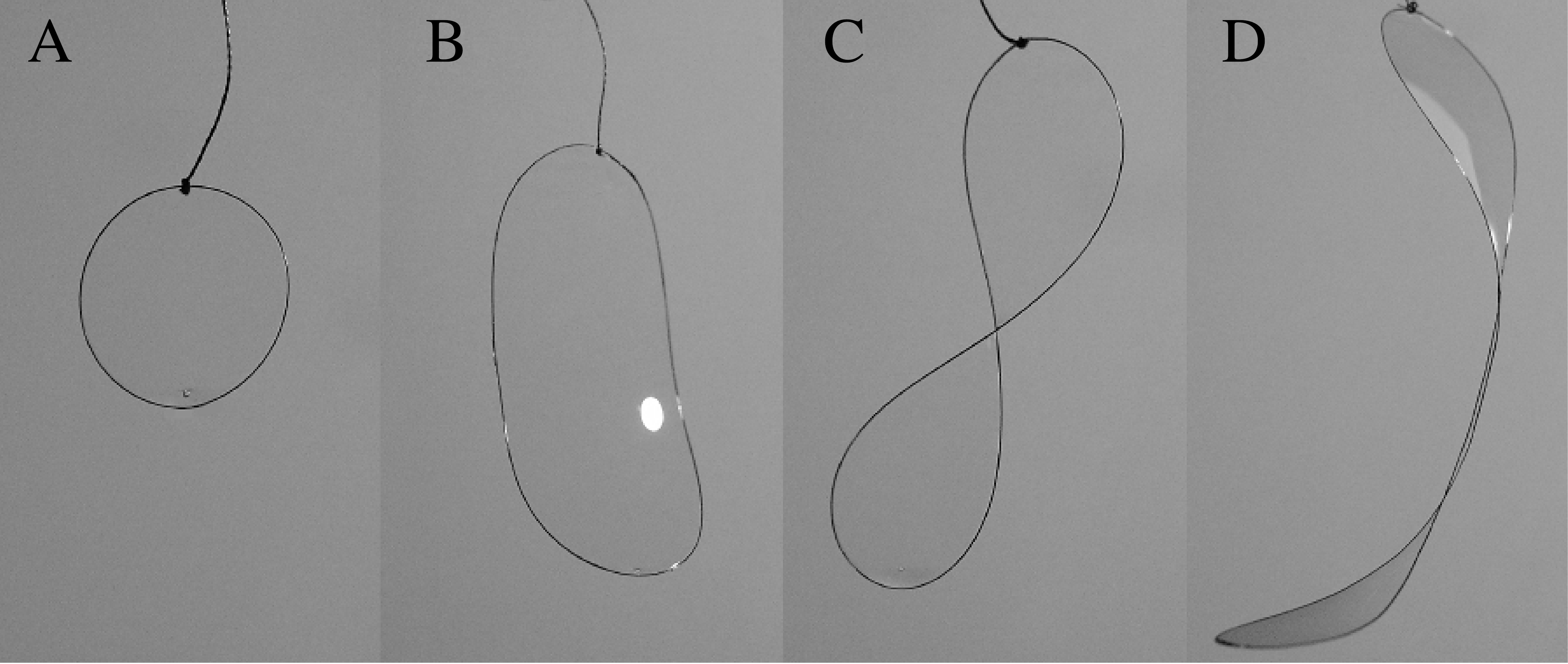}	
\caption{\label{fig:selection}A physical realization of the Euler-Plateau problem. (A)-(D) Examples of minimal surfaces bounded by an elastic line obtained by dipping (and removing) a loop of fishing line (0.3 mm diameter) in soapy water.  The interplay between the  forces due to surface tension of the film and the elastic forces arising at the boundary can be quantified via the dimensionless number $\gamma=\sigma L^{3}/\alpha$ with $\sigma$ the surface tension, $L$ the length and $\alpha$ the bending rigidity of the boundary.  For small values of $\gamma$, the soap film takes the form of (A) a planar disk, but upon increasing $\gamma$ it buckles in the form of (B) a saddle, (C) a twisted figure eight and (D) a two headed racket-like structure.}
\end{figure}

In Fig. \ref{fig:selection}A-D, we show a gallery of shapes  obtained via a simple experiment of dipping a loop of fishing-line into a solution of water and dish soap, and pulling it out. The resulting shapes are quite varied: short loops span a planar disk (Fig. \ref{fig:selection}A), intermediate loop lengths cause the film to spontaneously twist out of the plane (\ref{fig:selection}B) and eventually form a planar figure-eight shape (Fig. \ref{fig:selection}C). At the center of the eight (where the boundary crosses) the surface normal rotates by 180$^{\circ}$, and the spanning minimal surface is helicoid-like. Increasing the length even further leads to a a two-headed racket-like structure (Fig. \ref{fig:selection}D). 

To understand these shapes and the transitions between them, we start with a  mathematical formulation of the problem using a variational principle, given in section 2. We then briefly describe some general geometrical and physical aspects of the resulting Euler-Lagrange equations before describing the first symmetry-breaking bifurcation from a planar circular disc to a planar elliptical disc in section 3. In section 4 we describe numerical approximations to the secondary bifurcations from the planar disk to twisted strips that resemble our experiments. In section 5, we characterize these using an asymptotic theory using a minimal trial shape that captures the symmetries of the solutions seen in our computations, and finally conclude in section 6 with a discussion of some open problems that our study points towards.

\section{Mathematical Formulation}
To quantify the equilibrium shapes of these softly restrained soap films, we use a minimal mathematical model that poses the resulting forms as a results of the minimization of the following functional representing the total energy of the system:
\begin{equation}\label{eq:energy}
F[M] = \sigma \int_{M} dA\,+\oint_{\partial M}ds\,(\alpha\kappa^{2}+\beta)\,.
\end{equation}
Here $M$ is a mapping of a two-dimensional closed disk representing a film of surface tension $\sigma$ bounded by a flexible curve of circular cross-section with bending rigidity $2\alpha$. The boundary $\partial M$ is thus treated as an {\em Elastica} immersed in three-dimensional Euclidean space $\mathbb{R}^{3}$, and $\beta$ is a Lagrange multiplier, analogous to line tension, that enforces inextensibility of the boundary:
\begin{equation}\label{eq:inextensibility}
L = \oint_{\partial M}ds
\end{equation}
The assumptions that underlie our minimal model are as follows: (i)  the soap film itself has negligible thickness and is isotropic in its tangent plane, with an energy proportional to its area, (ii) the soap film meets the the bounding elastic filament at its edge along a line, so that it can apply normal forces to the filament, but no torques, (iii) the bounding elastic filament has a circular cross-section with a uniform bending and twisting stiffness along its entire length and its ends are glued together without twisting them relative to each other. 

Then, it follows that the filament must remain untwisted even when buckled and bent out of the plane, although it can and does have geometric torsion. This is because the translational symmetry of the circular cross-section along the filament implies that the twisting strain along it must be a constant (Love 1927); indeed the twist is the Noetherian invariant associated with this symmetry. Taken together with the fact that the closed filament has zero twist initially, this implies that the filament must always remain twist-free;  therefore the energetic contributions associated with twisting deformations vanish identically. As we will show, our simple theory suffices to capture the qualitative nature of the shapes and the transitions between them.

To derive the conditions for the extremum of the functional (\ref{eq:energy}), we will impose a small virtual deformation of the film and its boundary and use this to calculate the corresponding change in energy and the necessary conditions for the resulting shape to be a local minimum.  For the boundary curve and  changes therein, we use a parametrization in terms of the Frenet frame $\{\bm{t},\,\bm{n},\,\bm{b}\}$ of the tangent,  normal and binormal vector at each point of the curve (Kreyszig 1991). This orthogonal triad evolves in space according to the Frenet-Serret equations: 
\begin{equation}\label{eq:frenet-serret}
\bm{t}' = \kappa\bm{n}\,, \qquad 
\bm{n}' =-\kappa\bm{t}+\tau\bm{b}\,, \qquad
\bm{b}' =-\tau\bm{n}\,,
\end{equation}
where $\kappa$ and $\tau$ are respectively the curvature and torsion of the boundary curve and $(.)' = d(.)/ds$ in Eq. \eqref{eq:frenet-serret}  where $s$ is the arc-length. A linear variation of the boundary energy due to a small deformation can be easily calculated and leads to the following form:
\begin{equation}\label{eq:boundary-variation}
\delta \oint_{\partial M}ds\,(\alpha\kappa^{2}+\beta) 
= \oint_{\partial M} ds\,[\alpha(2\kappa''+\kappa^{3}-2\tau^{2}\kappa)-\beta\kappa]\xi + \alpha\oint_{\partial M}ds\,(4\kappa'\tau+2\kappa\tau')\eta
\end{equation}
where $\xi$ and $\eta$ are respectively the displacement along the normal and binormal vectors (see Mumford 1993 for a pedagogical derivation). The linear variation of the area bounded by the curve can be expressed in a standard coordinate basis $\bm{g}_{i}=\partial_{i}\bm{R}$ (with $\bm{R}$ position vector of the film and $\partial_{i}$ partial derivative with respect to the $i$-th coordinate) as (Lenz \& Lipowsky 2000):
\begin{equation}
\delta \int_{M} dA = \int_{M} dA\,(\nabla\cdot\bm{u}-2Hw) = \oint_{\partial M} ds\,\bm{\nu}\cdot\bm{u}-2\int_{M}dA\,H w\,.\label{eq:bulk-variation}
\end{equation}
Here $\bm{u}$ and $w$ are respectively the displacement along the tangent plane of the film (expressed in the basis $\bm{g}_{i}$) and its normal direction $\bm{N}$, $H$ is the mean curvature of the film and $\bm{\nu}=\bm{t}\times\bm{N}$ is the outward directed tangent vector normal to the boundary curve.

Compatibility demands that the bulk and boundary variations must be consistent with each other. To understand the conditions for this, we now need to express the first term on the right-hand side of Eq. \eqref{eq:bulk-variation} in the Frenet frame of the boundary curve. This can can be done by letting:
\begin{equation}
\bm{N} = \cos\vartheta\,\bm{n}+\sin\vartheta\,\bm{b}\,,
\end{equation}
where $\vartheta$ is {the} angle between the surface normal to the film and the normal to the boundary curve and will be herein referred to as {\em contact angle}. Standard algebraic manipulations then yields:
\begin{equation}\label{eq:contact-angle}
\bm{\nu}\cdot\bm{u} = \eta \cos\vartheta\,\bm{b}-\xi\sin\vartheta\,\bm{n}\,.
\end{equation}
Enforcing the condition that the energy variations \eqref{eq:boundary-variation} and \eqref{eq:bulk-variation} vanish while satisfying the compatibility condition \eqref{eq:contact-angle} yields:
\begin{subequations}\label{eq:euler-lagrange1}
\begin{gather}
\kappa''+\frac{1}{2}\kappa^{3}-\left(\tau^{2}+\frac{\beta}{2\alpha}\right)\kappa-\frac{\sigma}{2\alpha}\,\sin\vartheta = 0\,, \\[10pt]	
2\kappa'\tau+\kappa\tau'+\frac{\sigma}{2\alpha}\,\cos\vartheta = 0\,, \\[10pt]
H=0.
\end{gather}
\end{subequations}
Eqs. \eqref{eq:euler-lagrange1} can be made dimensionless by rescaling the arc-lenght $s$, the curvature $\kappa$ and the torsion $\tau$ by the length $L$ of the boundary. This identifies unambiguously a single dimensionless parameter $\gamma = \sigma L^{3}/\alpha$ representing the ratio between surface tension and bending rigidity. In sections 2 and 3 we will see how tuning $\gamma$ (for instance by changing the length of the boundary) triggers a sequence of buckling transitions, giving access to a variety of shapes. Determining the shape of the soap film then corresponds to solving (\ref{eq:euler-lagrange1}) subject to periodicity of the boundary curve; here, we note that it is not sufficient to impose periodicity on the curvature and torsion, but we must impose it directly on the coordinate description of the curve.

\section{Some general considerations and bifurcation from the flat state}

Before we proceed to understand the nature of the solutions of the Euler-Lagrange equations derived in the previous section, we provide some general considerations of the nature of the equations. From a geometrical perspective, since the projection of the surface normal $\bm{N}$ on the boundary normal and binormal direction are associated with the normal and geodesic curvature of the boundary, namely:
$\kappa_{n}=\kappa\cos\vartheta$ and $\kappa_{g}=\kappa\sin\vartheta$, an alternative way to write the equation satisfied by the boundary curve  \eqref{eq:euler-lagrange1} in terms of the {\em extrinsic} and {\em intrinsic} curvature of the boundary leads to
\begin{subequations}\label{eq:euler-lagrange2}
\begin{gather}
\kappa''+\frac{1}{2}\kappa^{3}-\left(\tau^{2}+\frac{\beta}{2\alpha}\right)\kappa-\frac{\sigma}{2\alpha}\,\frac{\kappa_{g}}{\kappa} = 0\,, \\[5pt] 2\kappa'\tau+\kappa\tau'+\frac{\sigma}{2\alpha}\,\frac{\kappa_{n}}{\kappa}= 0.
\end{gather}
\end{subequations}

From a physical perspective,  is interesting to note that Eqs. \eqref{eq:euler-lagrange1} can be also obtained from the classical equilibrium equations of rods in presence of a body force directed along $\bm{\nu}$. These require the total force and torque acting on a line element of the rod to vanish at equilibrium and yield (Landau \& Lifshitz 1986):
\begin{equation}\label{eq:mech-equilibrium}
\bm{F}'= - \bm{K}\,, \qquad\qquad\bm{M}'= \bm{F}\times\bm{t}\,.
\end{equation}
The vectors $\bm{F}$ and $\bm{M}$ can be expressed in the Frenet frame of the bent rod as:
\begin{subequations}\label{eq:force-torque}
\begin{gather}
\bm{F} = (\alpha\kappa^{2}-\beta)\,\bm{t}+2\alpha\kappa'\,\bm{n}+2\alpha\kappa\tau\,\bm{b}\,,\\[10pt]
\bm{M} = -2 \alpha\kappa\,\bm{b}\,.
\end{gather}
\end{subequations}
We see that while torque balance is automatically satisfied, taking the external force $\bm{K}=\sigma\bm{\nu}$ and differentiating the expression for $\bm{F}$ readily yields Eqs. \eqref{eq:euler-lagrange1}.

We now turn to a number of properties of the system that can be inferred from Eqs. \eqref{eq:euler-lagrange1} and \eqref{eq:euler-lagrange2} without knowledge of the exact solution. Multiplying Eq. (\ref{eq:euler-lagrange2}b) by $\kappa$ allows us to recast  the resulting expression as $(\kappa^{2}\tau)'=-(\alpha/2\sigma)\kappa_{n}$, which integrated along the length of the closed boundary  yields the following integral formula for the normal curvature $\kappa_{n}$:
\begin{equation}\label{eq:total-normal}
\oint ds\,\kappa_{n} = 0\,,
\end{equation}
valid for all values of the physical parameters. 

Furthermore, it is possible to consider a special class of solutions of Eqs. \eqref{eq:euler-lagrange1} for which the contact angle $\vartheta$ is constant along the curve. Then Eq. \eqref{eq:total-normal} can be used to prove that $\vartheta$ must then be $90^{\circ}$. To understand this latter statement we first need to recall that the normal vector $\bm{n}$ is undefined at inflection points (i.e. where $\kappa=0$). From Eq. \eqref{eq:contact-angle} it follows that, in order for $\vartheta$ to be defined everywhere along the boundary,  the curvature of the boundary must be non-zero everywhere, i.e. $\kappa > 0$. If $\vartheta$ is constant, writing $\kappa_{n}=\kappa\cos\vartheta$ in Eq. \eqref{eq:total-normal} yields:
\begin{equation}
\cos\vartheta\oint ds\,\kappa = 0\,, 
\end{equation}
but since $\kappa$ is strictly positive, this last condition is possible if and only if $\cos\vartheta=0$, thus $\vartheta=\pm \pi/2$: i.e. $\kappa_{n}=0$ and the surface normal must be perpendicular to the curve normal everywhere if $\vartheta$ is to be constant. This imposes strong constraints on the nature of the solution in this case as we will now argue.

Curves having zero normal curvature everywhere are called {\em asymptotic curves} and their tangent vectors {\em asymptotic directions}. Every curve that lies on the plane is clearly asymptotic. On a minimal surface, on the other hand, there are two asymptotic curves passing through each point where the Gaussian curvature $K$ is non-zero and these curves always intersect at 90$^{\circ}$. If the minimal surface contains flat points (i.e. isolated points where $K=0$), these are traversed by $n>2$ asymptotic curves crossing at $\pi/n$ and forming, in their neighborhood, $n$ valleys separated by ridges (Koch \& Fisher 1990). These observations suggest that the only simply connected minimal disk bounded by a closed asymptotic curve lies in the plane;  otherwise the asymptotic directions in the disk would form a vector field tangent to the boundary. In a simply connected disk this vector must then enclose a singularity or vortex, in the neighborhood of which the structure of the asymptotic directions would disagree with the previous classification. Hence, we conclude that the only solution having constant contact angle is given by the flat disk.

However, even in this case, it is possible for the soap film to be non-circular when the parameter $\gamma$ is large enough. To see this, we note that a trivial solution of Eqs. \eqref{eq:euler-lagrange1} is a circular disk whose radius $R$ satisfies the equation:
\begin{equation}\label{eq:disk}
\sigma R^{3}+\beta R^{2}-\alpha = 0\,,
\end{equation}
for all parameter values, with $\beta$ enforcing the condition $L=2\pi R$ with $L$ the length of the boundary. For large enough surface tension, however, we might expect this configuration to become unstable and the disk to buckle into a more complex shape. To test this hypothesis we analyze the stability of the disk with respect to a small periodic displacement in the radial direction $\delta R = \rho_{k}\sin k\phi$, where $\phi$ is the polar angle and $\rho_{k}$ a small amplitude. Expanding the energy \eqref{eq:energy} to second order in $\rho_{k}$ yields:
\begin{equation}\label{eq:energy-approx}
F \approx \sigma (\pi R^{2}) + 2\pi R \left(\frac{\alpha}{R^{2}}+\beta\right)+ \frac{1}{2}\pi\rho_{k}^{2}\left[\sigma+\frac{k^{2}\beta}{R}+\frac{\alpha(2k^{4}-5k^{2}+2)}{R^{3}}\right]\,.
\end{equation}
Here the Lagrange multiplier $\beta$ is given by Eq. \eqref{eq:disk} with $R=L/2\pi$, so that: $\beta = [\alpha-\sigma\,(L/2\pi)^{3}]/(L/2\pi)^{2}$. Introducing the dimensionless parameter $\gamma=\sigma L^{3}/\alpha$ and replacing $\beta$ in Eq. \eqref{eq:energy-approx}, it is easy to verify that the coefficient of the second order term becomes negative when:
\begin{equation}\label{eq:buckling}
\gamma > 16\pi^{3}(k^{2}-1)\,.
\end{equation}
Thus the first mode that goes unstable is that associated with an elliptical deformation with $k=2$, corresponding to the critical value $\gamma^{*} = 48\pi^{3}$. This elastic instability was first discussed by Love (1927) in the case of a circular ring constrained to lie on the plane and subject to a uniform normal pressure, and has since been the subject of rigorous bifurcation analysis (Flaherty {\em et al.} 1972). As $\gamma$ is increased still further, the non-circular shape eventually twists out of the plane. To understand these transitions, we resort to a combination of numerical and asymptotic   approximations.

\section{Numerical simulations of the soap film shapes}	

To explore the variety of possible shapes resulting from the solution of Eqs. \eqref{eq:euler-lagrange1} and the transitions between them, we  minimized a discrete {analog} of the total energy \eqref{eq:energy}. The soap film is approximated as simplicial complex consisting of an unstructured triangular mesh. The internal edges of the triangles are treated as elastic springs of zero rest-length so that the total energy of the mesh is given by:
\begin{equation}\label{eq:discrete-energy}
F[M] = \alpha\sum_{v \in \partial M} \langle s_{v} \rangle \kappa_{v}^{2} + k \sum_{e \in M } |e|^{2}\,.	
\end{equation}	
The first sum runs over the boundary vertices and $\langle s_{v} \rangle = (s_{v}+s_{v-1})/2$ is the average of the length of the two edges meeting at $v$. The curvature of the boundary is calculated as $\kappa_{v}=|\bm{t}_{v}-\bm{t}_{v-1}|/\langle s_{v} \rangle$ with $\bm{t}_{v}$ and $\bm{t}_{v-1}$ the tangent vectors at $v$. The second sum in Eq. \eqref{eq:discrete-energy} runs over all internal edges. If the triangles are equilateral, this yields a discrete approximation for the soap-film energy with the spring stiffness proportional to the surface tension, i.e. $k=\sqrt{3}\sigma/4$. {The choice of minimizing the squared length of the edges, instead of the area of the triangles, is motived by numerical stability. Replacing the spring energy in \eqref{eq:discrete-energy} with the sum of the area of the triangles, has the effect of shortening the range of interaction between the boundary and the interior to the single strip of triangles at the boundary. Most standard local optimization algorithms would then attempt to reduce the area of these triangles to zero, thus suppressing the interaction between the boundary and the interior.  In the continuum limit, the energy \eqref{eq:discrete-energy} approaches that given by \eqref{eq:energy} for  original problem, and  the sequence of shapes obtained with this method is in excellent agreement with our experimental observations shown in Fig. \ref{fig:selection}.

Numerical minimization of Eq. \eqref{eq:discrete-energy} using a conjugate gradient method leads to a variety of shapes depending on the values of $\alpha$, $k$, the boundary length $L$ as well as the initial shape of the domain. An interactive gallery of the shapes obtained based on our numerical simulations is available on-line\footnote{\url{http://www.seas.harvard.edu/softmat/Euler-Plateau-problem/}} and we invite the reader to explore these shapes. A first set of simulations were run by starting with an unstructured triangular mesh of 474 vertices bounded by an {\em elongated} hexagonal perimeter; the mesh points are started out randomly displaced transverse to the plane. As the spring constant $k$ was increased, corresponding to an increase in the  surface tension, the system moved through a series of transitions shown in Fig. \ref{fig:twisting}. Beyond a critical value of the surface tension (relative to the bending stiffness, as discussed earlier), the planar disk  buckled into a two-fold  mode consistent with our linear stability analysis. The geometry of this shape has been rigorously described (see for example {Flaherty {\em et al}. 1972}, Arreaga {\em et al}. 2002, Vassilev {\em et al}. 2009) and is related with problem of the collapse of a infinitely long cylindrical pipe under a uniform pressure difference. For $kL^{3}/\alpha\approx 643$,  the planar elliptical shape transitions to a twisted non-planar saddle-like shape, as experimentally observed. This is accompanied by an increase in the normal curvature of the boundary and there is  a progressive twisting of the central ``waist'' of the surface. A further increase in $k$ leads to an increase in the twist  until the surface normal undergoes a full 180$^{\circ}$ rotation across the waist for $kL^{3}/\alpha\approx 740$.  The surface then reverts to a planar configuration in the shape of the number eight (see Fig. \ref{fig:twisting}E). 

\begin{figure}[h!]	
\centering
\includegraphics[width=1\textwidth]{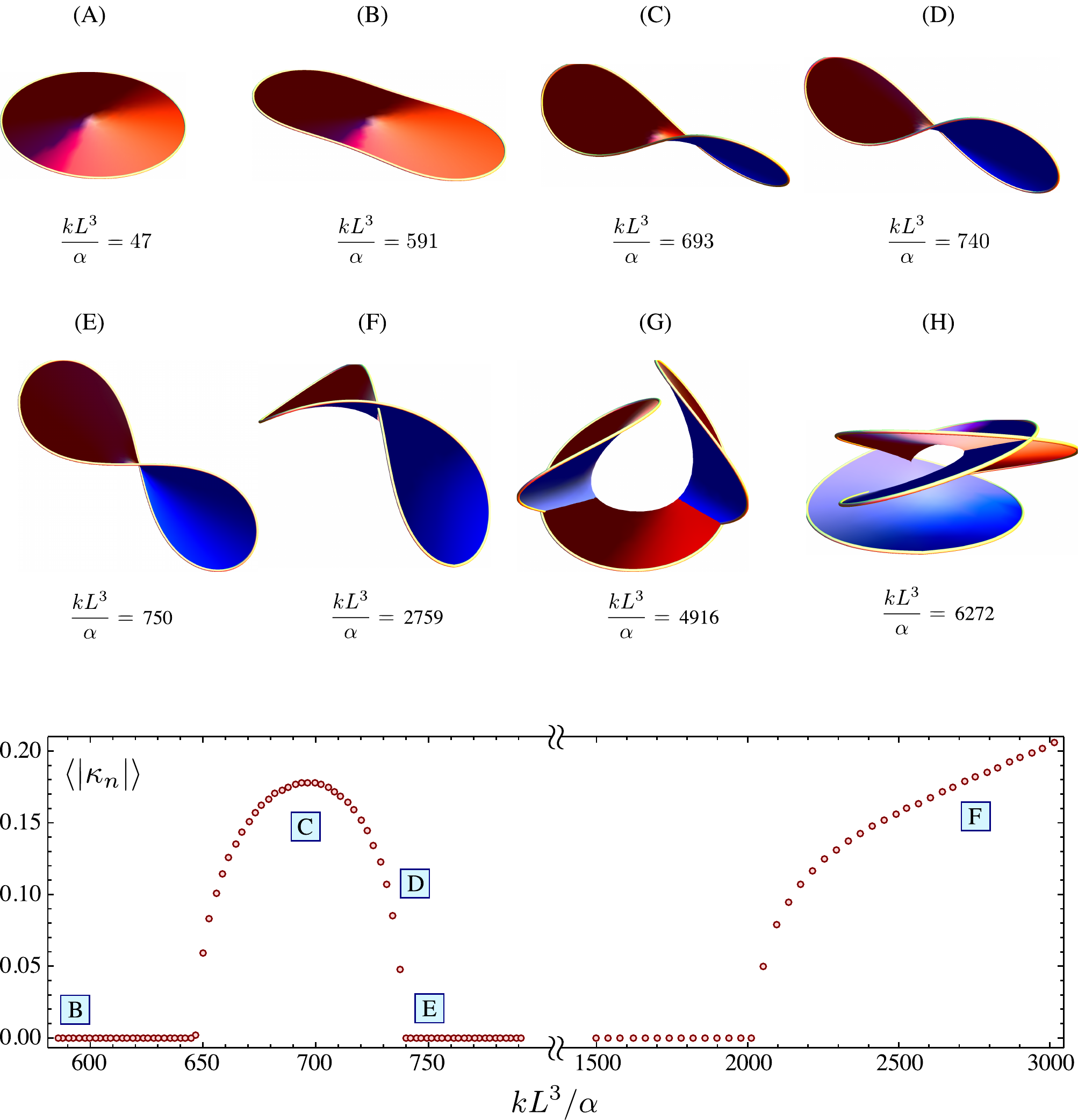}	
\caption{\label{fig:twisting}Some solutions of the Euler-Plateau problem and a bifurcation diagram. Shapes obtained by minimizing the discrete approximation \eqref{eq:discrete-energy} of the continuum energy of a minimal surface bounded by an elastic line \eqref{eq:energy} as a function of the dimensionless number $kL^{3}/\alpha$. For small values of the parameter, the flat circular disk of radius $L/2\pi$ has the lowest energy. Upon increasing $kL^{3}/\alpha$ the disk first buckles into flat elliptical shape (see text). For $kL^{3}/\alpha>643$, this elliptical disk  deforms into a three-dimensional saddle-like shape, becoming increasingly twisted  with an increase in  $kL^{3}/\alpha>643$ until the surface eventually becomes flat again, but now forming a fully twisted eight-shaped conformation. This twisted flat shape is a local energetic minimum for $740<kL^{3}/\alpha<2000$. For still larger $kL^{3}/\alpha$ the twisted surface self-intersects at the waist of the figure-eight and the two lateral lobes bend toward each other, while for very large $kL^{3}/\alpha$ the surface exhibits several self-intersections leading to the complex structure (H). Below the shapes we show the bifurcation diagram for the system, characterizing the shape using the absolute normal curvature of the boundary integrated along its length (a measure of the amplitude of the instability) as a function of the bifurcation parameter $kL^{3}/\alpha$, along with the location of the transitions described above.}
\end{figure}

Once the surface reaches the planar figure-eight conformation an increase in the surface tension does not produce further conformational changes until $kL^{3}/\alpha\approx 2000$. Beyond this value, the surface again becomes non-planar, self-intersects at the pinch of the eight, and the two lateral lobes start bending toward each other. The curvature of the lobes increases with increasing $k$ until  they intersect to produce the complex shape shown Fig. \ref{fig:twisting}H. The transition to this latter surface is analogous to the transition from the figure-eight to the two headed racket-like structure observed in our experimental realization of the problem. However, in the physical experiment, self-intersection cannot occur and adhesion favors the formation of a line of contact between the two lobes (see Fig. \ref{fig:selection}D). 

Not surprisingly, given the nonlinear nature of the governing equations (\ref{eq:euler-lagrange1}), the solutions are not expected to be unique, as our stability analysis of the flat disk has already indicated. In our numerical simulations, we also find that the final state is a strong function of the initial configuration of the system.  If the  initial configuration is {\em not elongated}, as in say a triangular mesh bounded by a hexagonal perimeter, which is closer to being circularly symmetric, we find other possible stable configurations. As in the previous case, for small values of the surface tension, the system rapidly relaxes into a flat circular disk. However, for $kL^{3}/\alpha\approx 855$ the system transitions to a saddle-like configuration that is the classical Enneper minimal surface. Upon increasing the surface tension still further, the curvature of the boundary and film becomes larger still and the surfaces self-intersect leading to the beautiful strucure shown in Fig. \ref{fig:enneper}D.

\begin{figure}[t]	
\centering
\includegraphics[width=1\textwidth]{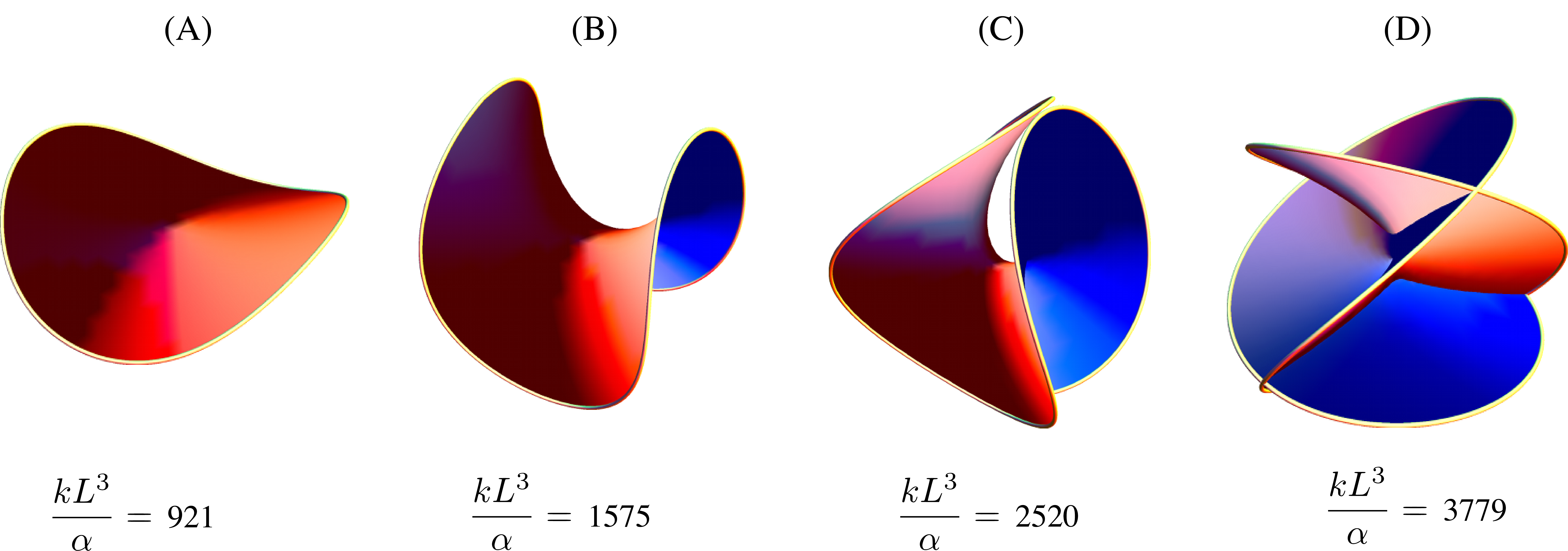}	
\caption{\label{fig:enneper}Further equilibrium solutions of the Euler-Plateau problem. An initial configuration of a triangular mesh bounded by a hexagonal perimeter gives the classic Enneper minimal surface rather than the figure-eight shape shown in Fig. 1.}
\end{figure}

\section{Asymptotic analysis}

\begin{figure}[t]
\centering
\includegraphics[width=1\textwidth]{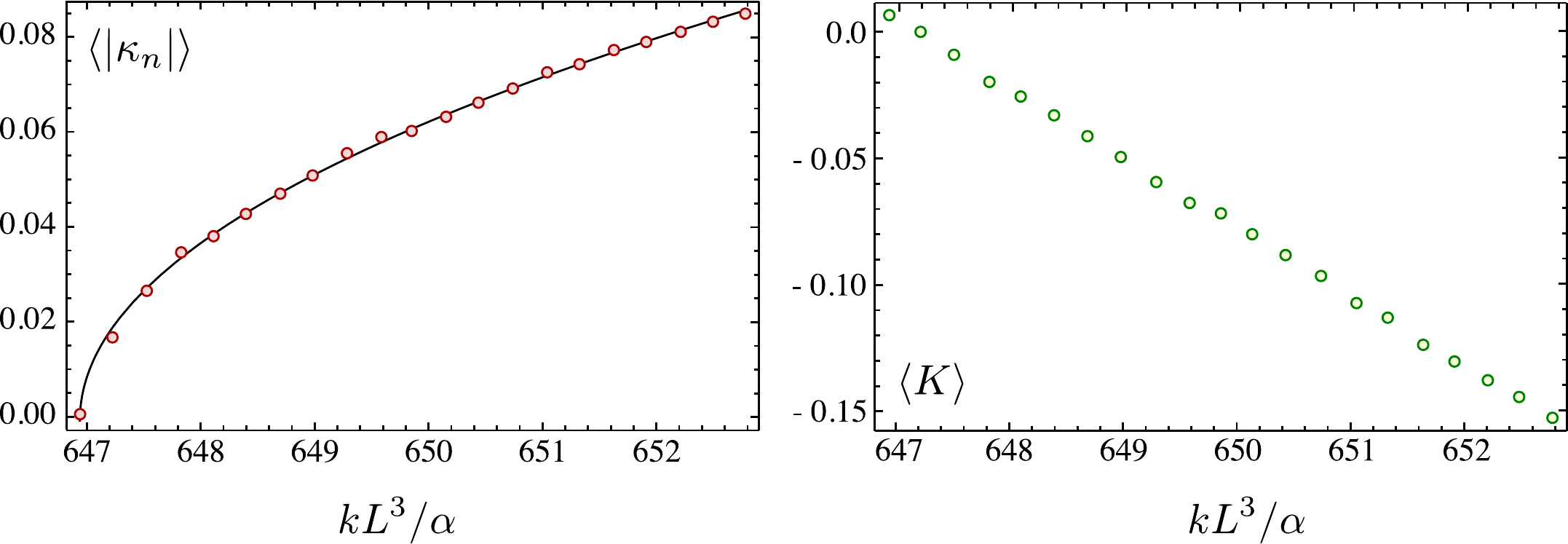}	
\caption{\label{fig:bifurcation}Characterizing the onset of twist and non-planarity. The left panel shows the average absolute normal curvature of the boundary $\langle |\kappa_{n}| \rangle$  as a function of $kL^{3}/\alpha$ in the neighborhood of the onset of the transition to non-planarity. The solid line on the left figure is obtained from a best fit of the data. The exponent resulting from the fit is $0.489 \pm 0.005$, and is consistent with the instability being of a supercritical pitchfork type. The right panel shows the average Gaussian curvature of the surface $\langle K \rangle$ as a function of $kL^{3}/\alpha$ in the neighborhood of the onset of the transition to non-planarity. The linear behavior is also consistent with our approximate analysis which yields the result shown in  the text.}
\end{figure}

The numerical data show that the transition from the planar two-fold symmetric shape to the non-planar twisted eight is consistent with a supercritical pitchfork bifurcation, with $\langle |\kappa_{n}|\rangle \sim (\gamma-\gamma^{*})^{1/2}$, where the angular brackets stand for an average along the boundary curve. Indeed Fig. \ref{fig:bifurcation} shows a best fit of the numerical data at the onset of the transition, the exponent obtained from the fit is $0.489 \pm 0.005$. 

To {understand} the origin of the pitchfork bifurcation, {it is helpful to introduce a specific approximate representation of the twisted soap film in terms of the following one-parameter family of surfaces}: 
\begin{equation}\label{eq:eight}
x = r\,(1+t^{2})\cos\phi\,,\quad
y = r\,(1-t^{2})\sin\phi\,,\quad
z = t\left(\frac{r^{2}}{R}\right)\sin 2\phi\,, 	
\end{equation}
where $r\in[0,R]$ and $\phi\in[0,2\pi]$ are the usual plane-polar coordinates. Here $t\in[-1,1]$ is a parameter that characterizes the family of shapes. {Eq. \eqref{eq:eight} is the simplest parameterization of a twisted saddle with disk topology and whose degree of non-planarity is controlled by the parameter $t$. For $t=0$, the surface given by (\ref{eq:eight}) is a disk of radius $R$; as $t$ becomes positive (negative), the surface becomes a right-handed (left-handed) twisted figure-eight, and for $t=\pm 1$ the surface reduces to a flat twisted figure-eight bounded by a lemniscate of Gerono. In the following we will exploit the topological equivalence between this family of surfaces and the actual twisted soap film to study the global properties of the surface at the onset of the transition, when $|t|\ll 1$.

The first fundamental form for the surface (\ref{eq:eight}) is given by:
\begin{subequations}\label{eq:metric-tensor}
\begin{gather}
g_{rr} = t^{4}+2t^{2}\left[\left(\frac{r}{R}\right)^{2}+\cos2\phi-\left(\frac{r}{R}\right)^{2}\cos4\phi\right]+1\\[5pt]
g_{r\phi} = 2rt^{2}\left[\left(\frac{r}{R}\right)^{2}\sin4\phi-\sin2\phi\right] \\[5pt]
g_{\phi\phi}/r^{2} = t^{4}+2t^{2}\left[\left(\frac{r}{R}\right)^{2}-\cos2\phi+\left(\frac{r}{R}\right)^{2}\cos4\phi\right]+1
\end{gather}
\end{subequations}
while its Gaussian curvature is given by:
\begin{equation}\label{eq:gaussian-curvature}
K = -\left(\frac{2t}{R}\right)^{2}+o(t^{3})
\end{equation}
For the boundary curve, its arc-length and squared curvature are given respectively by:
\begin{subequations}\label{eq:asymptotic-curve}
\begin{gather}
ds^{2} = R^{2}\left[(1+t^{2})^{2}-2t^{2}(\cos 2\phi - \cos 4\phi)\right] \\[10pt]
\kappa^{2} = \frac{1+t^{2}(10-6\cos 4\phi)-t^{4}(2-6\cos2\phi-2\cos6\phi)+t^{6}(10-6\cos 4\phi)+t^{8}}{R^{2}\left[(1+t^{2})^{2}-2t^{2}(\cos 2\phi - \cos 4\phi)\right]^{3}}
\end{gather}
\end{subequations}
Knowing the properties of the parametrized surface and its boundary in terms of its metric Eqs. \eqref{eq:metric-tensor} and curvature \eqref{eq:asymptotic-curve} we can approximate the total energy of the system given by (\ref{eq:energy}) as:
\begin{equation}\label{eq:approx-energy}
F \approx \frac{\pi \alpha}{R}\left[2+\frac{\sigma R^{3}}{\alpha}+t^{2}\left(10+\frac{\sigma R^{3}}{\alpha}\right)-t^{4}\left(9+\frac{5}{3}\frac{\sigma R^{3}}{\alpha}\right)\right]
\end{equation}}
Finally, the condition of inextensibility for the boundary curve is given by:
\begin{equation}\label{eq:length}
L = \oint ds = \pi R \left(2+2t^{2}-t^{4}\right)+o(t^{6})\,.
\end{equation}
Using Eq. (\ref{eq:length}) in Eq. \eqref{eq:approx-energy} to eliminate $R$ and then taking the derivative of the approximate energy Eq. \eqref{eq:approx-energy} with respect to $t$ yields the following equation of equilibrium for the amplitude parameter $t$ that characterizes the shape of the boundary:
\begin{equation}
t\left(96\pi^{3}-\gamma\right)+\frac{2}{3}\,\gamma t^{3}= 0\,,
\end{equation}
We recognize this as the normal form of a supercritical pitchfork bifurcation and see that it implies:
\begin{equation}\label{eq:t-critical}
t \sim (\gamma-\gamma^{*})^{1/2}\,,\qquad\gamma^{*} = 96\pi^{3}\,.
\end{equation}
We note that the critical value of the bifurcation parameter $\gamma^*$ is consistent with our exact linear stability analysis, but our approximate parametrization of the shape allows us to go beyond this via a weakly nonlinear analysis.

The normal and geodesic curvature of the boundary are given by:
\begin{subequations}\label{eq:kn-kg}
\begin{gather*}
\kappa_{n} = -\frac{2t}{R}\sin 2\phi+ \frac{t^{3}}{R}\,(3\sin 2\phi-\sin 4\phi+\sin 6\phi) + o(t^{4})\,, \\[7pt]
\kappa_{g} = \frac{1}{R}+\frac{t^2 (1+3 \cos 2\phi -5\cos 4\phi)}{2 R}+o(t^{4})\,,
\end{gather*}
\end{subequations}
which at leading order in $t$ leads to:
\begin{equation}\label{eq:kn-critical}
\langle |\kappa_{n}| \rangle \approx 8t \sim (\gamma-\gamma^{*})^{1/2}\,,
\end{equation}
and confirmed by our numerical simulations as shown in the left panel of Fig. \ref{fig:bifurcation}.

Having characterized the curvature of the elastic boundary curve, we now turn to the bounded minimal surface. The Gaussian curvature $K$ of the surface in the neighborhood of the onset of the bifurcation to non-planarity is given by Eq. \eqref{eq:gaussian-curvature},
\begin{equation}\label{eq:gauss-bonnet} 
\oint ds\,\kappa_{g}+\int dA\,K = 2\pi\,.
\end{equation}
The integrated geodesic curvature can be easily calculated at leading order in $t$:
\[
\oint ds\,\kappa_{g} \approx 2\pi(1+2t^{2})\,,
\]
from which, on using \eqref{eq:gauss-bonnet} and \eqref{eq:kn-critical}, we find that  the integrated Gaussian curvature $\langle K \rangle =\int dA\,K$ is :
\begin{equation}
\langle K \rangle \sim -(\gamma-\gamma^{*})\,,
\end{equation}
in excellent agreement with the numerical simulations shown in the bottom panel of Fig. \ref{fig:bifurcation}.

\section{Conclusions}

Through simple physical experiments using loops of flexible thread that are deformed by a soap film that is bounded by the loops, we have demonstrated a natural link between two of the oldest problems in geometrical physics and the calculus of variations, bringing together Euler's elastica and Plateau's problem.  The union of these leads to the question of the shape of  minimal surfaces bounded by elastic lines, a new class of questions at the nexus of geometry and physics. We formulate a simple geometric variational principle for the shapes of these surfaces and their boundaries and a corresponding set of coupled system of partial differential equations for the double free boundary problem.  A combination of scaling concepts, asymptotic analysis and numerical simulation allows us to characterize the qualitative nature of the solutions that result.

However, many questions remain. While a minimal physical realization of the Euler-Plateau problem is a kitchen-sink experiment, its ramifications are likely to go far beyond this specific  manifestation, just as the humble soap film and the elastic filament have been relevant for the study of matter not just at the everyday scale, but also for systems that range from molecules to black holes (Thomas {\em et al}. 1988, Kamien 2002, Penrose 1973), while inspiring art (Ferguson). In addition, both Euler's elastica and Plateau's problem have inspired a variety of  mathematical developments in topology, geometry, analysis and beyond (Bryant \& Griffiths 1983, Langer \& Singer 1984, Osserman 2002, Colding \& Minicozzi 2007, Morgan 2008). Perhaps their fruitful union might also do the same ? Some immediate mathematical questions include proofs of existence, regularity of solutions as well as measures of the non-uniqueness of solutions particularly for large values of the only parameter in the problem $\gamma$ or its discrete analog $kL^3/\alpha$. In terms of applications, the Euler-Plateau problem is likely lurking behind almost every physical manifestation of the Plateau problem, itself an idealization; time alone will uncover these. Further generalizations that introduce the effects of finite shearing and bending rigidity in the soap film, the torques from the soap film's Plateau borders on the filament, the effects of twist deformations of the filament due to anisotropy and inhomogeneity of its cross-section etc. are just some avenues that might be worth exploring; indeed some of these arise naturally in problems of vertebrate gut morphogenesis (Savin et al. 2011) and elsewhere. 

\section*{Acknowledgments}

We acknowledge support from the NSF Harvard MRSEC, the Harvard Kavli Institute for Nanobio Science \& Technology, the Wyss Institute and the MacArthur Foundation. We are grateful to Aryesh Mukherjee for help with the experiment and William Meeks III for discussions.

\end{document}